\documentclass{amsart}
\usepackage{amssymb}
\usepackage{amsmath}
\usepackage{amsthm}
\usepackage{mathtools}
\usepackage{ifpdf}
\usepackage{graphicx}

\pdfoutput=1
\ifpdf
\usepackage[pdftex]{hyperref}
\else
\usepackage[hypertex]{hyperref}
\fi

\hypersetup{
citecolor=green,
menucolor=red,
citebordercolor=0 1 0,
linkbordercolor=1 0 0,
pdfpagemode=UseOutlines,
pdftitle={Kink-soliton trains of nonlinear Schroedinger equations},
pdfauthor={Stefan Le Coz <slecoz@math.univ-toulouse.fr>, Tai-Peng Tsai <ttsai@math.ubc.ca>},
pdfsubject={Kink-soliton trains of nonlinear Schroedinger equations},
pdfkeywords={soliton, multi-soliton, kink-soliton train, infinite soliton train, nonlinear Schroedinger equations} 
}

\newtheorem{theorem}{Theorem}

\newtheorem{assumption}{Assumption}

\DeclarePairedDelimiter{\norm}{\lVert}{\rVert}

\newcommand{\psld}[2]{\left( #1,#2 \right)_{2}}

\newcommand{\dual}[2]{\left\langle #1,#2 \right\rangle}

\newcommand{\eps}{\varepsilon}

\newcommand{\N}{\mathbb{N}}
\newcommand{\R}{\mathbb{R}}

\newcommand{\hu}{H^1(\R^d)}

\renewcommand{\leq}{\leqslant}
\renewcommand{\geq}{\geqslant}

\begin{document}

\title[Kink-soliton trains of nonlinear Schr\"odinger equations]{Finite and infinite soliton and kink-soliton
trains of nonlinear Schr\"odinger equations}

\author[S.~Le Coz]{Stefan Le Coz}
\author[T.-P.~Tsai]{Tai-Peng Tsai}

\address{Institut de Math\'ematiques de Toulouse,
\newline\indent
Universit\'e Paul Sabatier
\newline\indent
118 route de Narbonne, 31062 Toulouse Cedex 9
\newline\indent
France}
\email{slecoz@math.univ-toulouse.fr}

\address{
Department of Mathematics,
\newline\indent
University of British Columbia,
\newline\indent
Vancouver BC
\newline\indent
Canada V6T 1Z2}
\email{ttsai@math.ubc.ca}

\thanks{
The research of S. Le Coz is supported in part by the french ANR
through project ESONSE.  The research of Tsai is supported in part by
NSERC grant 261356-13 (Canada).
}

\date{\today}

\maketitle

\begin{abstract}
  We will first review known results on multi-solitons of dispersive partial differential equations, which are special solutions behaving like the sum of many weakly-interacting solitary waves. We will then describe our recent joint work with Dong Li on nonlinear Schr\"odinger equations: Assuming the composing solitons have sufficiently large relative speeds, we prove the existence and uniqueness of a soliton train which is a multi-soliton composed of infinitely many solitons. In the 1D case, we can add to the infinite train an additional half-kink, which is a solution with a non-zero background at minus infinity. 
\end{abstract}

\section{Introduction}\label{sec:introduction}

The nonlinear Schr\"odinger equation
\begin{equation}\label{eq:nls}\tag{{\sc nls}}
\left\{
\begin{aligned}
&iu_t+\Delta u+f(u)=0,\\
&u(t=0)=u_0,
\end{aligned}
\right.
\end{equation}
where $u$ is a complex valued function on $\R\times\R^d$ and the nonlinearity $f:\mathbb C\to\mathbb C$ is phase
covariant, i.e. there exists $g:\R^+\to\R$ such that $f(z)=g(|z|^2)z$, 
appears in various physical contexts, for example in nonlinear optics or in the modelling of Bose-Einstein condensates. Mathematically speaking, it is one of the model nonlinear dispersive PDE, along with the Korteweg-De Vries equation and the nonlinear wave equation. 

The local Cauchy theory in the energy space $H^1(\R^d)$ for \eqref{eq:nls} is well understood (see e.g. \cite{Ca03} and the references cited therein). In this survey we are interested in the long time dynamics of global solutions. Essentially, two effects may be at play at large time. First of all, if the nonlinearity is not too strong, the linear part of the equation can dominate and solutions may behave as if they were solutions to the free linear Schr\"odinger equation. This is the \emph{scattering effect}. On the other hand, in some cases the nonlinear term dominates and the solution tends to concentrate, with possible blow-up in finite time. This is the \emph{focusing effect}. At the equilibrium between these two effects, one may encounter many different types of structures that neither scatter nor focus. The most common of these non-scattering global structures are the solitons, but there exist also dark solitons, kinks, etc. A generic conjecture for nonlinear dispersive PDE is the \emph{Soliton Resolution Conjecture}. Roughly speaking, it says that, as can be observed in physical settings, any global solution will eventually decompose at large time into a scattering  part and well separated non-scattering structures, usually a sum of solitons. Apart from integrable cases (see e.g. \cite{ZaSh72}), such conjecture is usually out of reach (see nevertheless the recent breakthrough \cite{DuKeMe12} on energy critical wave equation). Intermediate steps toward  this conjecture are existence and stability results of configurations with well separated non-scattering structure, like multi-solitons, multi-kinks, infinite soliton and kink-soliton trains, etc. 
Our purpose in this paper is  to review some of the existing results on this topic.

\section{Solitons, dark solitons and kinks}\label{sec:2}
We review in this section some of the known elementary non-scattering structure solutions of \eqref{eq:nls}.

Let us first consider special solitons, the standing waves. A \emph{standing wave} is a solution to \eqref{eq:nls} of the form $u(t,x)=e^{i\omega t}\phi(x)$, where $\omega \in \R$ and $0\not=\phi\in\hu$ is a localized solution to the elliptic stationary equation
\begin{equation}\label{eq:snls}
-\Delta \phi+\omega\phi-f(\phi)=0.
\end{equation}
Since the ground work of Berestycki and Lions \cite{BeLi83-1,BeLi83-2}, it is well known that \eqref{eq:snls} admits solutions in $\hu$ if $\omega>0$ and $f(z)=g(|z|^2)z$ verify the following hypotheses.

\begin{assumption}\label{as:1}
(energy-subcritical) Let $d\ge 1$. Suppose $f(z)=g(|z|^2)z$ with $g \in C^0([0,\infty), \R) \cap C^2((0,\infty),\R)$, $g(0)=0$, 
\begin{align*} 
 |sg^{\prime}(s)| +|s^2 g^{\prime\prime}(s)| \le C \cdot (s^{\alpha_1/2} +s^{\alpha_2/2}), \qquad\forall\, s>0,
\end{align*}
where $C>0$, $0<\alpha_1 \le \alpha_2< {\alpha_{\max}}$, $\alpha_{\max}=+\infty$ if $d=1,2,$ $\alpha_{\max}=\frac{4}{d-2}$ if $d\geq 3$.\\
 (focusing) There exists $s_0>0$, such that
\begin{equation*}
G(s_0) := \int_0^{s_0} g(\tilde s) d\tilde s > \omega s_0. 
\end{equation*}
\end{assumption}

The profile solutions of \eqref{eq:snls} in $\hu$ are in general called \emph{bound states}. Among bound states, it is common to distinguish between the \emph{ground states} and the \emph{excited states}. Recall first that three quantities are conserved along the $H^1$-flow of \eqref{eq:nls}: the energy, the mass and the momentum, defined as follows.
\begin{gather*}
E(u):=\frac{1}{2}\norm{\nabla u}_2^2-\int_{\R^d}F(u)dx,\qquad F(z):=\int_0^{|z|}f(s)ds,\\
M(u):=\frac{1}{2}\norm{u}_2^2,\qquad
P(u):=\Im\int_{\R^d}\bar{u}\nabla udx.
\end{gather*}
The ground states are minimizers of the action related to \eqref{eq:snls} (defined by $S=E+\omega M$) and are in general positive, radial and unique (see  \cite{GiNiNi81,Kw89} and the recent progresses  \cite{ByJeMa09,Ma09}). For $d=1$, there exist only ground states, whereas for $d\geq 2$ there are infinitely many excited states.

A \emph{soliton} is a standing wave of \eqref{eq:nls} that has been given a speed thanks to a Galilean transform. Since \eqref{eq:nls} is Galilean invariant, a soliton is still a solution of \eqref{eq:nls}. Explicitly, a soliton with frequency $\omega>0$, speed $v\in\R^d$, initial phase $\gamma \in \R$  and position $x_0\in\R^d$ has the form 
\[
R_{\phi, \omega, \gamma, x_0,v}:= \phi(x-vt-x_0 ) \exp
\left( i \Bigl( \frac 12 v \cdot x - \frac 14 |v|^2 t + \omega
t +\gamma \Bigr) \right).
\]
The dynamical properties of solitons are mostly known when the nonlinearity $f$ is of power-type $f(z)=|z|^\alpha z$ and the profile is a ground state (we will refer to such solitons as \emph{ground state solitons}, opposed to \emph{excited state soliton} or in general \emph{bound state soliton}). For $L^2$-subcritical $\alpha$ ($ \alpha<\frac 4d$), the ground state soliton is orbitally stable \cite{CaLi82} (i.e. stable up to phase shifts and translations), whereas for $L^2$-critical and supercritical $\alpha$ ($\alpha=\frac4d$, $\alpha>\frac4d$) it is unstable by blow-up \cite{BeCa81,We83} (i.e. there exists an initial data in any neighborhood of the soliton such that the corresponding solution to \eqref{eq:nls} blows up in finite time). In general, excited state solitons are expected to be unstable \cite{Gr88,Jo88,Mi05-1,Mi05-2,Mi07}.

For the focusing power-type nonlinearities, solitons are the only non-scattering solutions with a fixed profile, which is always localized. However, for other nonlinearities,  there exist also other types of non-scattering solutions with a fixed profile, for example the so-called dark and grey solitons or the kinks. 

A \emph{dark soliton} is a travelling wave solution of \eqref{eq:nls} of  the form $\phi(x-ct)$, where $c\in\R^d$ is the velocity and $\phi$ is a profile which has the particularity to be non-localized but with a constant modulus at infinity. The analysis of such type of solitons is much less developped than for localized (\emph{bright}) solitons. Most of the works deal with the case where \eqref{eq:nls} is the Gross-Pitaevskii equation, i.e. when $f(z)=(1-|z|^2)z$. In this case, existence of dark solitons in dimension $1$ follows from direct computations. One may refer to \cite{Ch12,Ch13} and the reference cited therein for a study of the existence and stability of dark solitons with generic nonlinearities in dimension $1$.  In higher dimension, existence of dark solitons has been a long time open problem. In the recent breakthrough \cite{Ma13}, it has been proved that, for generic nonlinearities and in dimension $d\geq 3$, dark solitons exist for any speed between $0$ and the speed of sound ($\sqrt{2}$ in the case of the Gross-Pitaevskii nonlinearity). 

A \emph{kink} is a soliton-type solution to \eqref{eq:nls} when $d=1$, 
\[
K_{\phi_K, \omega, c,\gamma, x_0,v}:= \phi_K\left(x-(c+v)t-x_0 \right) \exp
\left( i \Bigl( \frac 12 v \cdot x - \frac 14 |v|^2 t + \omega
t +\gamma \Bigr) \right),
\]
 but with a profile $\phi_K$ which has different limits at $-\infty$ and $+\infty$. Here $c$ is the intrinsic velocity associated to $\phi_K$.
 When $f$ is the Gross-Pitaevskii nonlinearity, i.e. $f(z)=(1-|z|^2)z$, there exists an explicit family of kink solutions (which are particular cases of dark solitons) given by 
$K(t,x)=\phi_{K}(x-ct)$,
\begin{equation}
\label{GP-kink} \phi_K(x)=\sqrt{\frac{2-c^2}{2}}\operatorname{tanh}\left(\frac{x\sqrt{2-c^2}}{2}\right)+i\frac{c}{\sqrt{2}},\quad 
 |c|<\sqrt{2}.
\end{equation}

In this paper, we will be particularly interested in cases where the kink profile $\phi$ has different limits \emph{with different modulus} at $-\infty$ and $+\infty$. This will be the case if the nonlinearity $f$ verifies the following assumption.
\begin{assumption}\label{as:kink}
For some $\omega_0>0$, there is a first $b>0$ such that
for $h(s)=\omega_0 s - f(s)$,
\begin{equation*}
h(b)=0,\quad h'(b)>0, \quad \int_0^b h(s)ds =0.
\end{equation*}
\end{assumption}
Under Assumption \ref{as:kink}, there exists a (unique up to translation) kink-profile $\phi_K\in \mathcal C^2(\R)$ with zero intrinsic velocity such that
\begin{equation*}
\left\{
\begin{aligned}
&-\phi_K''+\omega_0\phi_K-f(\phi_K)=0,\\
&\lim_{x\to -\infty}\phi(x)=b,\quad \lim_{x\to+\infty}\phi(x)=0.
\end{aligned}
\right.
\end{equation*}
To our knowledge, except in our works \cite{LeLiTs13,LeTs13}, such kinks have never been investigated in the analysis of the long time behavior of solutions to \eqref{eq:nls}. One reason for that is that they are genuinely infinite energy solutions and that no proper renormalization exists to make them energy finite. 

\section{Multi-solitons, infinite soliton trains and soliton-kink solutions}

As mentioned in Section \ref{sec:introduction}, it is expected that global solutions to \eqref{eq:nls} will eventually decompose into a scattering part and well separated non-scattering structures. In this section, we review the existence and stability results for solutions composed of several of the elementary non-scattering structures described in Section \ref{sec:2}. 
The basic example is the so-called \emph{multi-soliton}, a solution of \eqref{eq:nls} build upon a finite number of solitons, 
\begin{figure}
\centering
\includegraphics{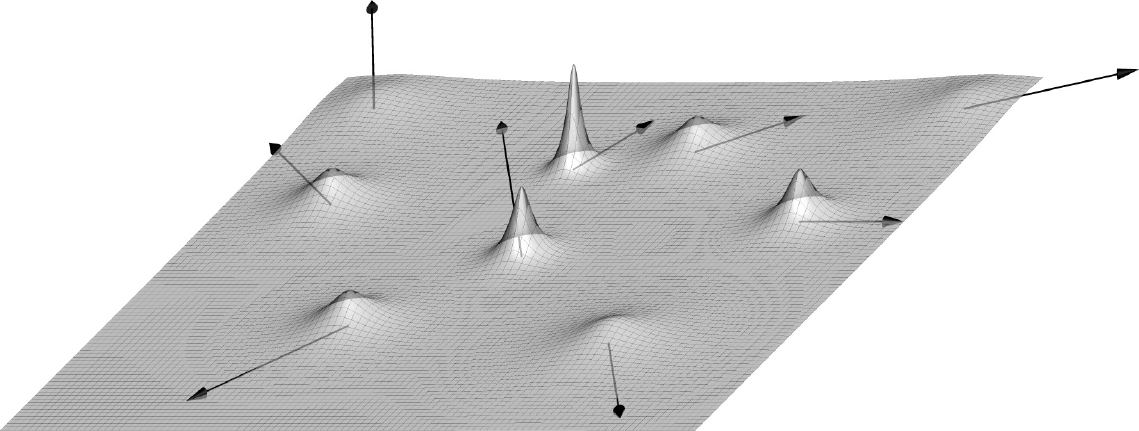}
\caption{Schematic representation of a multi-soliton}
\label{fig:multisolitons}
\end{figure}
Figure \ref{fig:multisolitons} contains a schematic representation of a multi-soliton.
To be a little more precise, let 
\begin{equation}\label{eq:profile}
R(t,x)= \sum_{j=1}^N R_{\phi_j,\omega_j,\gamma_j,x_j,v_j}(t,x) =:
\sum_{j=1}^N R_j(t,x),
\end{equation}
where each $R_j$ is a soliton made from some parameters $(\omega_j,\gamma_j, x_j, v_j)$ and bound state $\phi_j$. As \eqref{eq:nls} is a nonlinear equation, $R$ is not a solution. What we call a \emph{multi-soliton} is a solution $u$ of \eqref{eq:nls} such that
\begin{equation}\label{eq:multi-soliton}
\lim_{t\to+\infty}\norm{u-R}_{X([t,+\infty)\times\R^d)}=0,
\end{equation}
where $X$ is some space-time function space, e.g. $L^\infty([0,+\infty),L^2(\R^d))$.

\subsection{The integrable cases}

The first result of existence of multi-solitons was obtained in Zakharov and Shabat \cite{ZaSh72} in the case of the $1$-d focusing cubic  (i.e. $d=1$, $f(z)=|z|^2z$) nonlinear Schr\"odinger equation via the inverse scattering method. Indeed, in this particular case the equation is completely integrable and one can obtain multi-solitons in a rather explicit manner. The companion situation to the $1$-d cubic case is the Gross-Pitaevskii case, i.e. $d=1$ and $f(z)=(1-|z|^2)z$. In that case \eqref{eq:nls} is again completely integrable, however no localized solitons exist. The counter parts of the solitons in this case are the dark solitons \eqref{GP-kink}, which have modulus $1$ at infinity. In this situation, it is also possible to show via the inverse scattering transform the existence of a solution of \eqref{eq:nls} behaving at large time as decoupled well-separated dark solitons. Note that in that case, due to the non-zero condition at infinity, the profile cannot be given simply as the sum of the dark solitons.

 For the $1$-d cubic \eqref{eq:nls}, Kamvissis \cite{Ka95} showed that it is possible  to push the inverse scattering analysis forward and obtain the existence of an\emph{ infinite soliton train}, i.e. a solution $u$ of \eqref{eq:nls} defined as in \eqref{eq:multi-soliton} but with $N=+\infty$ in the definition of $R$. In fact, the result given in \cite{Ka95} is sharper: it is shown that, if $d=1$, $f(z)=|z|^2z$ and under some technical hypotheses,  any solution to \eqref{eq:nls} with initial data in the Schwartz class will eventually decompose at large time as an infinite soliton train and a ``background radiation component''. The existence of an infinite train of dark solitons has not been proved  for the Gross-Pitaevskii case. 

Other exotic solutions are available for $1$-d cubic NLS, for example {\it breather} solutions 
which are periodic in space and whose magnitudes approach uniform values as $ |t| \to \infty$, see e.g.~\cite{MR915545}. 

\subsection{Existence of multi-solitons, the energy method}

Apart from the case previously mentioned, \eqref{eq:nls} is not completely integrable and it is generically not possible to apply the inverse scattering method to obtain the existence of multi-solitons. The first existence result of multi-solitons in a non-integrable setting was obtained by Merle in \cite{Me90} as a by-product of the proof of existence of multiple blow-up points solutions for $L^2$-critical \eqref{eq:nls}, i.e. $f(z)=|z|^\frac4dz$. The techniques initiated in \cite{Me90} were then developed in \cite{CoLC11,CoMaMe11,MaMe06} to obtain the following result. 

\begin{theorem}[Existence of multi-soliton by energy method, \cite{CoLC11,CoMaMe11,MaMe06}]\label{thm:energy}
Assume $f(z)=|z|^\alpha z$ with $0<\alpha<\alpha_{\max}$. 
Let $R$ be the profile given in \eqref{eq:profile}.
Let $\omega_\star$ and $v_\star$ be given by 
\[
\omega_\star=\frac{1}{2}\min\left\{ \omega_j,j=1,...,N \right\},\quad
v_\star = \inf_{j,k=1,\ldots,N, j\not = k}  |v_j -v_k | 
.
\]
There exist $\mu=\mu(d,N)$ and 
$v_\sharp:=v_\sharp(\phi_1,...,\phi_N)\geq0$ such that if
$v_\star>v_\sharp$ then there exist $T_0\in\R$ and a solution 
$u\in \mathcal{C}([T_0,+\infty),\hu)$ of \eqref{eq:nls} satisfying
\[ 
\norm{u(t)-R(t)}_{H^1}\leq e^{-\mu\sqrt{\omega_\star} v_\star t}, \quad \forall t\ge
T_0.
\]
If in addition all $\phi_j$ are ground states, then the result holds with $v_\sharp=0$,
i.e., any $v_*>0$.
\end{theorem}

For the sake of simplicity, we have stated the result only for power-type nonlinearities, but its first part is in fact valid for any $\mathcal C^1$ nonlinearity verifying hypotheses a little stronger than Assumption \ref{as:1}. As mentioned before, the case $\alpha=\frac4d$ was treated by Merle \cite{Me90}, the ground state case by Martel and Merle \cite{MaMe06} for $\alpha<\frac4d$ and C\^ote, Martel and Merle \cite{CoMaMe11} for $\alpha>\frac4d$, and the excited state case by C\^ote and Le Coz \cite{CoLC11}. 

The proofs in \cite{CoLC11,CoMaMe11,MaMe06} follow a similar scheme. The idea is to choose an increasing sequence of time $(T^n)$ with $T^n\uparrow+\infty$ and consider the solutions $(u_n)$ to \eqref{eq:nls} which solve the equation backward in time with final data $u_n(T^n)=R(T^n)$. The sequence $(u_n)$ is an approximate sequence for a multi-soliton. To show its convergence, two arguments are at play. First, one shows that there exists a time $T_0$ independent of $n$ such that $u_n$ satisfies on $[T_0,T^n]$ the \emph{uniform estimates}
\[
\norm{(u_n-R)(t)}_{H^1}\leq e^{-\mu\sqrt{\omega_\star}v_\star t}.
\]
Second, we have \emph{compactness} of the sequence of initial data $u_n(T_0)$, i.e. there exists $u_0\in\hu$ such 
\[
u_n\to u_0\quad\text{strongly in }H^s(\R^d), \quad 0\leq s<1.
\]

To prove the uniform estimates, one first recall that the second derivative of the action around a soliton is coercive up to some $L^2$ scalar products. Precisely, given a soliton $R_0$ with parameters $(\omega_0,\gamma_0,x_0,v_0)$, the action is given by
\[
S=E+\left(\omega_0+\frac{|v_0|^2}{4}\right)M+v_0\cdot P
\]
and there exist $K\in\N$ and $(\xi^0_k)_{k=1,\ldots,K}$ such that we have for any $\eps\in\hu$
\[
\dual{S''(R_0)\eps}{\eps}\gtrsim\norm{\eps}_{H^1}^2-\sum_{k=1}^K
|\psld{\eps}{\xi_k^0}|^2.
\]
Hence to control the difference $\eps=u_n-R$, it is enough to construct a functional $\mathcal S$ similar to  $S$ but adapted to  the multisoliton profile $R$ and to get rid of the bad $L^2$ scalar products. The construction of the functional is done by gluing together each functional $S_j$ suitably localized. The localization is possible since each bound state, henceforth each soliton, is exponentially decaying at infinity. 

Getting rid of the bad $L^2$ scalar products is the trickiest part. In \cite{MaMe06}, the authors modulated the solitons in scaling, translation and phase to cancel the scalar products. This was possible only because they were in the $L^2$-subcritical case with ground states. In \cite{CoMaMe11}, the authors used a topological argument to select a final data which was not $R(T^n)$ but close enough to it, and such that for this final data the bad scalar products  were vanishing. In \cite{CoLC11}, a bootstrap argument on the $L^2$-norm of $\eps$ was used to get an a priori control of the type $|\psld{\eps}{\xi_k^0}|\leq \frac{1}{v_\star}\norm{\eps}_2$, hence allowing to control the scalar products for $v_\star$ large enough. 

The compactness argument is based on the virial identity and does not present major difficulty. 

The energy method is very flexible and can be adapted to many other situations, e.g. for multi-solitons of  Klein-Gordon equations \cite{BeGhLe13,CoMu12} or multi-speeds solitary waves of Schr\"odinger systems \cite{IaLe12}. It suits very well situations with a finite number of well localized solitons with finite energy. However, its implementation is far from being trivial when the number of solitons is infinite or when one soliton is replaced by a kink. In \cite{LeLiTs13,LeTs13}, the authors have developped an  approach suitable to situations where the energy technique fails to be directly applicable.

\subsection{Existence of infinite trains and kink-solitons solutions, the fixed point argument}

To prove the existence of solutions which have a priori infinite energy such as multi-solitons with kinks attached at both ends or infinite trains of solitons, an approach based on a fixed point argument around the desired profile has been followed in \cite{LeLiTs13,LeTs13}. 

We briefly recall the definition of the Strichartz space.
For dimension $d\ge 1$,  a pair of exponents $(q,r)$ is said
\emph{(Schr\"odinger) admissible} if
\begin{equation}
\frac 2 q + \frac d r = \frac d2, \qquad 2\le q,r \le \infty, \; \text{and } (d,q,r)\ne (2,2,\infty).
\end{equation}
On $I\times \mathbb R^d$ where $I\subset \R$, we define the Strichartz norm
\begin{equation} \label{strichartz_def}
\| u \|_{S(I)} :=  \sup_{\text{$(q,r)$ admissible}} \|u \|_{L_t^{q} L_x^r (I\times \mathbb R^d)}.
\end{equation}
For $d=2$, we also need to impose $q>q_1$ in the above norm for some $q_1$ slightly larger than $2$ to stay away from the forbidden endpoint. The Strichartz space  $S(I)$ is the closure of all test functions in $\R \times \R^d$ under this norm. We denote by
$N(I)$ the dual space of $S(I)$.

Let us first state a result on the existence of infinite soliton trains. 
Its conditions are far from optimal for simplicity of presentation.
The same method can be applied to construct finite and infinite soliton trains for more general nonlinearities.

\begin{theorem}[Existence of infinite soliton trains \cite{LeLiTs13,LeTs13}] \label{thm:infinite}
Let $d \ge 1$.
Assume $f(z)=|z|^\alpha z+g_2(|z|^2)z$ where 
$g_2 \in C^0([0,\infty), \R) \cap C^2((0,\infty),\R)$, $g(0)=0$, 
\begin{align*} 
 |sg_2^{\prime}(s)| +|s^2 g_2^{\prime\prime}(s)| \le C (s^{\alpha_{1.5}/2} +s^{\alpha_2/2}), \qquad\forall\, s>0,
\end{align*}
where $C>0$, $0<\alpha < \alpha_{1.5} \le \alpha_2< {\alpha_{\max}}$. There 
exist $r_0>\max(1,\frac{d\alpha}2) $, $c_1>0$ and $v_\sharp \gg 1$ such that, if an infinite soliton train profile
$R_\infty$ is given by 
\[
R_\infty = \sum_{j=1}^{\infty} R_j,
\]
with parameters $\omega_j>0$,
$\gamma_j\in \R$, $x_j=0$, $v_j\in \R^d$ satisfying
\begin{itemize}
\item (uniform bound for bound states) for some $0<a<1$ and $C>0$,
\begin{align*}
|\phi_j(x)|+ \omega_j^{-1/2} |\nabla \phi_j(x)| \le C \omega_j^{1/\alpha} e^{-a \omega_j^{1/2}|x|}, \quad \forall x\in\R^d, \forall j\in\mathbb N,
\end{align*}
\item (Integrability)
\begin{equation*}
\sum_{j=1}^\infty 
 \omega_j^{\frac 1{\alpha} -\frac d {2r_0}} <\infty,
 \end{equation*}
\item (High relative speeds) 
\begin{align*}
v_* = \inf_{j,k\in \N, j\not = k} \sqrt { \omega_j  } |v_k -v_j | 
\ge v_\sharp   ,  
\end{align*}
\item (Gradient bound)%
\begin{align*}
V_* = \sum_{j \in \N} \langle v_j \rangle  \omega_j^{\frac 1{\alpha} -\frac d {4}} <\infty\quad \text{if 
} \alpha < \frac {\alpha_2}{2+\alpha_2}.
\end{align*}
\end{itemize}
Then there exists  a
solution $u$ of \eqref{eq:nls} 
satisfying, for some $T_0\ge 0$,
\begin{equation*}%
\|u-R_\infty\|_{S([t,\infty))} \le e^{-c_1 v_\star  t}, \quad \forall t \ge T_0.
\end{equation*}
Moreover, such profiles $R_\infty$ do exist for every such nonlinearity $f$.
\end{theorem}

Note that, in the integrable case,  it was proved that an integrability condition on the parameters of the infinite trains is necessary for their existence, see \cite[Remark 2]{Ka95}.

\begin{figure}
\centering
\includegraphics{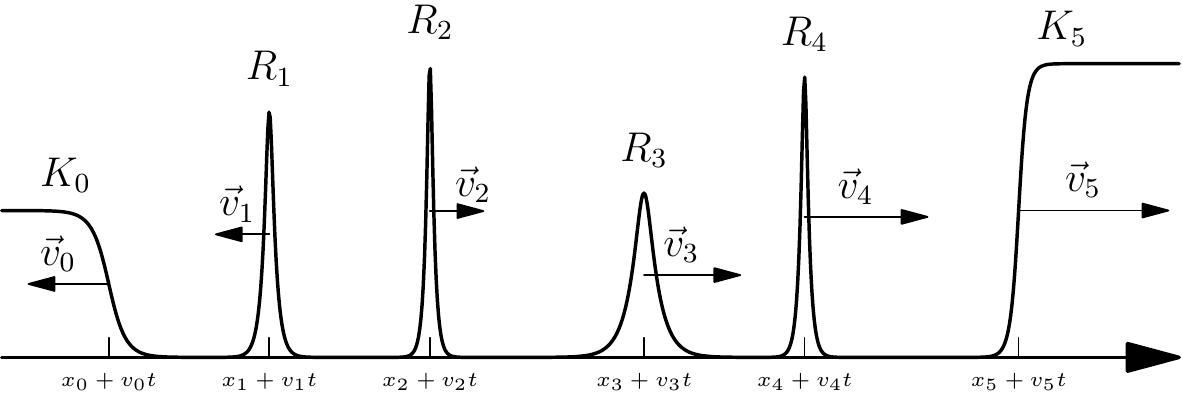}
\caption{Schematic representation of a kink-soliton train}
\label{fig:multi-kink}
\end{figure}

We now switch to dimension $d=1$ to investigate the existence of solutions of \eqref{eq:nls} composed of kinks and solitons as represented in Figure \ref{fig:multi-kink}.

\begin{theorem}[Existence of finite kink-soliton  trains \cite{LeLiTs13}] \label{thm:kink}
Let $d=1$. 
Assume that $f$ is an energy subcritical nonlinearity such that there exists two kink profiles $\phi_{K_0},\phi_{K_{N+1}}$ with $\phi_{K_0}'$ and $\phi_{K_{N+1}}'$ exponentially localized and the boundary conditions
\begin{align*}
& \lim_{x\to-\infty}\phi_{K_0}(x)\neq 0,&&	 \lim_{x\to+\infty}\phi_{K_0}(x)= 0,\\
& \lim_{x\to-\infty}\phi_{K_{N+1}}(x)= 0,&&	 \lim_{x\to+\infty}\phi_{K_{N+1}}(x)\neq 0.
\end{align*}
Assume also the existence of  $N$ soliton profiles $\phi_1,\dots,\phi_N$. Define the profile
\[
W(t,x):=K_0(t,x)+\sum_{j=1}^NR_j(t,x)+K_{N+1}(t,x).
\]
with parameters $(v_j,x_j,\omega_j,\gamma_j)_{j=0,\dots,N+1}\subset\R^4$ such that $v_0<\dots<v_{N+1}$.
Define $v_\star$ by 
\[
v_\star:=\inf\{|v_j-v_k|;\;j,k=0,\dots,N+1,\,j\neq k\}.
\]
Then there exist $v_\sharp>0$ (independent of $(v_j)$) large enough, $T_0\gg 1$
and constants $\mu_1,\mu_2>0$ such that if $v_\star>v_\sharp$, then
there exists a (unique) multi-kink solution
$u\in\mathcal{C}([T_0,+\infty),H^1_{\mathrm{loc}}(\mathbb R))$ to
  \eqref{eq:nls} satisfying on $[T_0,+\infty)$ the estimate
\[
e^{\mu_1v_\star t}\norm{u-W}_{S([t,+\infty))}+e^{\mu_2v_\star t}\norm{\nabla(u-W)}_{S([t,+\infty))}
\leq 1.
\]
\end{theorem}
In the above theorem, one kink can be dropped.

We may also have a kink attached to an infinite soliton train. For simplicity we choose a special nonlinearity. See \cite{LeTs13} for general assumptions on $f$.

\begin{theorem}[Existence of infinite kink-soliton trains \cite{LeTs13}] \label{thm:kink2}
Let $d=1$ and 
\[
f(u)=|u|^\alpha u - |u|^\beta u,\quad 0<\alpha<\beta<\infty.
\]
It satisfies Assumption \ref{as:kink} and there is a kink profile $\phi_K$.
If either $0<\alpha<4/3$, or $4/3\le \alpha < \sqrt 2<\beta=2/\alpha$, then 
there exist $r_0>1$, $c_1>0$ and $v_\sharp\gg 1$ such that, if an infinite kink-soliton train profile 
\[
W=K+R_\infty = K+\sum_{j=1}^{\infty}R_j
\]
has its parameters satisfying the same assumptions as in Theorem \ref{thm:infinite} (with index $j$ starting from $0$),
then there exists  a
solution $u$ of \eqref{eq:nls} 
satisfying, for some $T_0\ge 0$,
\begin{equation*}%
\|u-W\|_{S([t,\infty))} \le e^{-c_1 v_\star  t}, \quad \forall t \ge T_0.
\end{equation*}
Moreover, such profiles $W$ do exist.
\end{theorem}

The strategy for the proofs of Theorems \ref{thm:infinite}, \ref{thm:kink}  and \ref{thm:kink2} is the following. 
Let $W$ be a profile around which we want to build a solution. $W$ can be an infinite train, a kink-soliton train, etc.
Since $W$ may be badly localized, we look for a solution of \eqref{eq:nls} in the form $u= W+\eta$, where $\eta$
satisfies the perturbation equation
\begin{equation}\label{eq:perturbation}
i\partial_t \eta  + \Delta \eta  + f( W+\eta )-f(W)+H=0,
\end{equation}
where $H$ is a source term coming from the fact that $W$ is not an exact solution, e.g. $H=f(R_\infty)-\sum_{j=1}^{+\infty}f(R_j)$ in the case of an infinite soliton train.
In Duhamel formulation, the
perturbation equation for $\eta $ becomes
\begin{align*}
\eta (t)= -i \int_t^{\infty} e^{i(t-\tau)\Delta}\Bigl( f(W+\eta )
-f(W)+H \Bigr) d\tau, %
\end{align*}
and the core of the proof is to perform a fixed point argument for this formulation. Here, two approaches are possible:  One is based on 
a combination of the dispersive estimate
\begin{align*}
\| e^{it\Delta} u \|_{p} \lesssim |t|^{-d(\frac 12 -\frac 1p)} \| u\|_{\frac p{p-1}}, \qquad\forall\, t\ne 0,\quad  \forall  p\in[2,+ \infty],
\end{align*}
and Strichartz estimates 
\[
\|u \|_{S((t_0,\infty))} \lesssim \|u_0\|_2 + \| F \|_{N((t_0,\infty))}\quad\text{for}\quad i\partial_t u +\Delta u =F,\quad u(t_0)=u_0.
\]
The other uses only Strichartz estimates, but for both $u$ and $\nabla u$.

\subsection{Open problems}

We conclude this paper by reviewing some open problems related to the exotic solutions presented here. %

\subsubsection{Uniqueness}
In Theorem \ref{thm:energy}, no uniqueness is proved for the multi-soliton. 
In Theorems  \ref{thm:infinite}, \ref{thm:kink} and  \ref{thm:kink2}, the solution presented is unique \emph{in the class of solutions satisfying a strong decay estimate towards the desired profile}; it does not preclude the possibility of existence of the same type of solution, but with a weaker decay towards the desired profile. In fact, it was proved in \cite{CoLC11} that as soon as one of the composing soliton is linearly unstable, then there exists a one parameter family of multi-solitons converging toward the same profile. Hence in such cases uniqueness is not to be expected. However, one may hope that classification results hold, e.g. one may expect that for $L^2$-supercritical power-type nonlinearities, the multi-solitons converging toward a fixed $N$-sum of solitons form a $N$-parameter family. Such results were obtained for the Korteweg-de Vries equation by Combet \cite{Co11}. See also \cite{Co10} for partial results in that direction for \eqref{eq:nls}. 
For the kink-solitons solutions or the infinite trains, no result is available yet.

\subsubsection{Stability}

Another natural question coming to mind when investigating the exotic solutions of nonlinear Schr\"odinger equation is their stability.  Again, the only available results concern finite multisolitons and the problem is completely open for infinite trains or kink-soliton solutions. For power-type nonlinearities, the only case where orbital stability of a multi-soliton holds is the Gross-Pitaevskii case, as has been proved by B\'ethuel, Gravejat and Smets for multi-dark-solitons  in \cite{BeGrSm12}. Stability of multi-solitons has been proved under restrictive hypotheses in \cite{MaMeTs06} for orbital stability and in \cite{Pe97,Pe04,RoScSo03} for asymptotic stability. The hypotheses on the nonlinearity (e.g. high regularity or flatness at the origin) exclude in particular the power-type nonlinearities. 

\subsubsection{Multikinks}

One particular feature of the kinks considered in Theorem \ref{thm:kink} is that they  converge to $0$ on one side. There exist however nonlinearities such that the kinks are connecting non-zero constants on both sides, for example $-1$ and $1$ for the black solitons of Gross-Pitaevskii. There may also exist situations with kinks connecting e.g. $0$ to $1$ and kinks connecting $1$ to $2$, etc. One would expect that it is possible to construct solutions by gluing together those kinks to get a solution with  modulus having increasing terraces shape. Such solutions have never been constructed for nonlinear Schr\"odinger equation. One reason for that is the current lack of appropriate ansatz for  such terraces shape solutions. 

Many other related open problems exist, e.g. constructing solutions when one of the composing element is a line-soliton in $2$-d.

 \bibliographystyle{abbrv}
 \bibliography{biblio}

\end{document}